\documentclass[12pt]{amsart}
\usepackage{amssymb}
\newtheorem{theorem}{Theorem}[section]
\newtheorem{lemma}[theorem]{Lemma}
\newtheorem{proposition}[theorem]{Proposition}

\theoremstyle{definition}
\newtheorem{definition}[theorem]{Definition}

\theoremstyle{remark}

\numberwithin{equation}{section}
\begin{document}
\title {Orthogonal Constant Mappings in Isosceles Orthogonal Spaces}
\author{Madjid Mirzavaziri}
\address{M. Mirzavaziri: Department of Mathematics, Ferdowsi University, P. O. Box 1159, Mashhad 91775, Iran \newline 
and \newline Banach Mathematical Research Group (BMRG)}
\email{mirzavaziri@math.um.ac.ir}
\author{Mohammad Sal Moslehian}
\address{M. S. Moslehian: Department of Mathematics, Ferdowsi University, P. O. Box 1159, Mashhad 91775, Iran \newline and \newline Centre of Excellency in Analysis on
Algebraic Structures (CEAAS), Ferdowsi Univ., Iran.}
\email{moslehian@ferdowsi.um.ac.ir}
\keywords{}

\begin{abstract}
In this paper we introduce the notion of orthogonally constant
mapping in an isosceles orthogonal space and establish stability
of orthogonally constant mappings. As an application, we discuss
the orthogonal stability of the Pexiderized quadratic equation
$f(x+y)+g(x+y)=h(x)+k(y)$.
\end{abstract}
\maketitle

\section{Introduction}

We say a functional equation $(\mathcal E)$ is {\it stable} if any
function $g$ satisfying the equation $(\mathcal E)$
``approximately'' is near to an exact solution of $(\mathcal E)$.

The stability problem of functional equations originated from a
question of Ulam \cite{ULA} concerning algebra homomorphisms. In
1941, Hyers \cite{HYE} gave a partial affirmative answer to the
question of Ulam in the context of Banach spaces. In 1978, Th.M.
Rassias \cite{RAS1} extended the theorem of Hyers. The result of
Rassias has provided a lot of influence in the development of
what we now call {\it Hyers--Ulam--Rassias stability} of
functional equations. During the last decades several stability
problems of functional equations have been investigated in the
spirit of Hyers--Ulam--Rassias. The reader is referred to
\cite{B-M1, CZE, H-I-R, JUN, MOS1, MOS3, RAS2} and references
therein for more comprehensive information on stability of
functional equations.

There are several concepts of orthogonality such as
Birkhoff--James, Phythagorean, isosceles, Singer, Roberts,
Diminnie, Carlsson, $T$-orthogonality, R\" atz, etc, in an
arbitrary real normed space ${\mathcal X}$, which can be regarded
as generalizations of orthogonality in the inner product spaces,
in general. These are of intrinsic geometric interest and have
been studied by many mathematicians; see \cite{B-M2}. Among them we deal with
isosceles orthogonality $\perp$. This notion was introduced by C.
James \cite{JAM} as follows: $x\perp y$ if and only if $\|x+y\| =
\|x-y\|$. He proved that a Banach space ${\mathcal X}$ of three
or more dimensions is an inner-product space if isosceles
orthogonality in ${\mathcal X}$ is homogeneous in one variable,
or equivalently, it is additive for biorthogonal pairs of vectors
(see also \cite{ALO1, ALO2, D-F-A}) . In this paper, a real normed
space ${\mathcal X}$ endowed with the isosceles orthogonality
$\perp$ is called an isosceles orthogonal space.

The orthogonally quadratic equation $f(x+y)+f(x-y)=2f(x)+2f(y),~~
x\perp y$ was first studied by Vajzovi\' c \cite{VAJ} where $\perp$ means the Hilbert space
orthogonality. Later, Drljevi\' c \cite{DRL}, Fochi \cite{FOC}
and Szab\' o \cite{SZA} generalized this result. An investigation
of the orthogonal stability of the Pexiderized quadratic equation
$f(x+y) + f(x-y)=2g(x)+2h(y),~~x\perp y$ may be found in \cite{M-M, MOS2}.

In this paper we introduce the notion of orthogonally constant
mapping in an isosceles orthogonal space and establish stability
of orthogonally constant mappings. As an application, we use this
notion in a natural fashion to discuss the orthogonal stability of
the Pexiderized quadratic equation $f(x+y)+g(x+y)=h(x)+k(y)$.

Throughout this paper, $({\mathcal X},\perp)$ denotes an isosceles
ortohgonal space and ${\mathcal Y}$ is a normed space. We also
denote the even and odd parts of a given function $\rho$ by
$\rho^e(x):=\frac{\rho(x)+\rho(-x)}{2}$ and
$\rho^o(x):=\frac{\rho(x)-\rho(-x)}{2}$, respectively.

\section{Orthogonal Mappings}

We start our work with the following definition.

\begin{definition}
(i) A mapping $c:{\mathcal X}\to {\mathcal Y}$ is called \textit{orthogonally
constant} if $c(x+y)=c(x-y)$ for all $x,y\in {\mathcal X}$ with $x\perp y$.

(ii) A mapping $f:{\mathcal X}\to {\mathcal Y}$ is called
\textit{approximately orthogonally constant} if there is a
positive number $\varepsilon > 0$ such that
$\|f(x+y)-f(x-y)\|\leq\varepsilon$ for all $x,y\in {\mathcal X}$
with $x\perp y$.
\end{definition}

\begin{lemma}\label{equinorm} Let $c:{\mathcal X}\to {\mathcal Y}$ be an orthogonally constant
mapping. If $x, y \in {\mathcal X}$ and $\|x\|=\|y\|$, then $c(x) = c(y)$.
\end{lemma}
\begin{proof} Let $\|x\| = \|y\|$. Setting $h=\frac{x+y}{2}$ and
$k=\frac{x-y}{2}$. Then $h\perp k$ and so $c(x)=c(h+k)=c(h-k)=c(y)$.
\end{proof}

\begin{proposition} Suppose that $c:{\mathcal X}\to {\mathcal Y}$ is an orthogonally constant
mapping. Then there is a mapping $g:{\mathbb R}\to {\mathcal Y}$ such that
$c(x)=g(\|x\|)$ for each $x\in{\mathcal X}$.
\end{proposition}
\begin{proof} Let $x_0\neq0$ be a fixed element of ${\mathcal X}$.
Define $g:{\mathbb R} \to {\mathcal Y}$ by $g(r) = c(\frac{r x_0}{\|x_0\|})$.
Then, by Lemma \ref{equinorm}, we have
\[g(\|x\|)=c(\frac{\|x\|x_0}{\|x_0\|})=c(x),\]
since $\|\frac{\|x\|x_0}{\|x_0\|}\|=\|x\|$.
\end{proof}

\begin{proposition}\label{approx}
Let $f:{\mathcal X}\to {\mathcal Y}$ be an approximately
orthogonally constant mapping such that
\[\|f(x+y)-f(x-y)\|\leq\varepsilon\]
for some $\varepsilon > 0$ and for all $x,y\in {\mathcal X}$ with
$x\perp y$. Then there is an orthogonally constant mapping
$c:{\mathcal X}\to {\mathcal Y}$ such that
\[\|f(x)-c(x)\|\leq \varepsilon,\]
for all $x\in{\mathcal X}$.
\end{proposition}
\begin{proof} Let $x_0\neq 0$ be a fixed element of ${\mathcal X}$.
Define $c:{\mathcal X}\to {\mathcal Y}$ by
$c(x)=f(\frac{x_0\|x\|}{\|x_0\|})$. Then $c$ is an
orthogonally constant mapping since $x \perp y$ implies that $\|x + eta \| = \|x - y\|$ and so
$f(\frac{x_0\|x + y\|}{\|x_0\|}) = f(\frac{x_0\|x - y\|}{\|x_0\|})$, whence
$c(x + y) = c(x - y)$. In addition
\[\|f(x)-c(x)\| = \|f(x)-f(\frac{x_0\|x\|}{\|x_0\|})\|\leq \varepsilon,\]
since
$\frac{1}{2}(x+\frac{x_0\|x\|}{\|x_0\|})\perp \frac{1}{2}(x-\frac{x_0\|x\|}{\|x_0\|})$.
\end{proof}

\section{Application to Pexiderized Quadratic Equation}

The problem of orthogonal stability of the quadratic equation
$f(x+y)+f(x+y)=2f(x)+2f(y)$ was discussed in \cite{MOS2}. Towards
a satisfactory study of the orthogonal stability of the
Pexiderized quadratic equation $f(x+y)+g(x+y)=h(x)+k(y)$ we use
the notion of approximately orthogonally constant in a natural
fashion as follows.

An approximately orthogonally Cauchy mapping is a mapping $f:
{\mathcal X} \to {\mathcal Y}$ for which there is a positive
number $\varepsilon
> 0$ such that
\begin{eqnarray*}
\|f(x+y)-f(x)-f(y)\|\leq \varepsilon,
\end{eqnarray*}
for all $x, y\in {\mathcal X}$ with $x\perp y$. By an
approximately orthogonally quadratic mapping we mean a mapping $f:
{\mathcal X} \to {\mathcal Y}$ for which there is a positive
number $\varepsilon > 0$ such that
\begin{eqnarray*}
\|f(x+y)+f(x-y)-2f(x)-2f(y)\|\leq \varepsilon,
\end{eqnarray*}
for all $x, y\in {\mathcal X}$ with $x\perp y$.

\begin{theorem} Suppose that ${\mathcal X}$ is a isosceles orthogonal space with
an orthogonal relation $\perp$ and ${\mathcal Y}$ is a Banach
space. Let the mappings $f, g, h, k:{\mathcal X}\to {\mathcal Y}$
satisfy $f(0) = g(0) = h(0) = k(0) = 0$ and the following inequality
\begin{eqnarray}\label{PQ}
\|f(x+y)+g(x-y)-h(x)-k(y)\|\leq\varepsilon,
\end{eqnarray}
for all $x, y\in {\mathcal X}$ with $x\perp y$. Then $f$ is a
linear combination of an approximately orthogonally Cauchy
mapping, an approximately orthogonally quadratic mapping and an
approximately orthogonally constant mapping. The same is true for
$g$, as well.
\end{theorem}
\begin{proof}
Set $\ell(x)=\frac{h(x)+k(x)}{2}$. If $x\perp y$ then $-x\perp
-y$, hence we can replace $x$ by $-x$ and $y$ by $-y$ in
(\ref{PQ}) to obtain
\begin{eqnarray}\label{-PQ}
\|f(-x-y)+g(-x+y)-h(-x)-k(-y)\|\leq \varepsilon.
\end{eqnarray}
By virtue of triangle inequality and (\ref{PQ}) and (\ref{-PQ}) we have
\begin{eqnarray}\label{o}
\|f^o(x+y)+g^o(x-y)-h^o(x)-k^o(y)\|\leq \varepsilon,
\end{eqnarray}
\begin{eqnarray}\label{e}
\|f^e(x+y)+g^e(x-y)-h^e(x)-k^e(y)\|\leq \varepsilon,
\end{eqnarray}
for all $x,y\in {\mathcal X}$.

Let $x\perp y$. Then $y\perp x$, and by (\ref{o})
\begin{eqnarray}\label{yperpx}
\|f^o(x+y)-g^o(x-y)-h^o(y)-k^o(x)\|\leq \varepsilon.
\end{eqnarray}
It follows from (\ref{o}) and (\ref{yperpx}) that
\begin{eqnarray}\label{fhkhk}
&&\|2f^o(x+y)-h^o(x)-k^o(x)-h^o(y)-k^o(y)\|\nonumber\\&\leq&\|f^o(x+y)+g^o(x-y)-h^o(x)-k^o(y)\|\nonumber\\
&&+\|f^o(x+y)-g^o(x-y)-h^o(y)-k^o(x)\|\nonumber\\
&\leq& 2\varepsilon.
\end{eqnarray}
for all  $x,y\in {\mathcal X}$ with $x\perp y$. In particular, for
arbitrary $x$ and $y=0$ we get
\begin{eqnarray}\label{fhk}
\|2f^o(x)-h^o(x)-k^o(x))\|&\leq& 2\varepsilon.
\end{eqnarray}
By (\ref{fhkhk}) and (\ref{fhk}), we have
\begin{eqnarray}\label{fo}
\|f^o(x+y)-f^o(x)-f^o(y)\|&\leq&\frac{1}{2}\|2f^o(x+y)-h^o(x)-k^o(x)-h^o(y)-k^o(y)\|\nonumber\\
&&+\frac{1}{2}\|2f^o(x)-h^o(x)-k^o(x))\|\nonumber\\
&&+\frac{1}{2}\|2f^o(y)-h^o(y)-k^o(y))\|\nonumber\\
&\leq& 3\varepsilon
\end{eqnarray}
for all  $x,y\in {\mathcal X}$ with $x\perp y$. Hence $f^o$ is
approximately orthogonally Cauchy mapping. It follows from
(\ref{fo}) and (\ref{fhk}) that $\ell^0$ is an approximately
orthogonally Cauchy mapping. Since $x \perp y$ implies that $x
\perp -y$, it follows from (\ref{o}) that
\begin{eqnarray*}
\|g^o(x+y)+ f^o(x+y)-h^o(x)- (-k)^o(y)\|\leq \varepsilon,
\end{eqnarray*}
hence by the same reasoning as above, we conclude that $g^o$ is
also approximately orthogonally Cauchy mapping.

Now, putting $x=0$ in \ref{e} we get
\[\|f^e(y)+g^e(-y)-k^e(y)\|\leq\varepsilon,\]
and putting $y=0$ in \ref{e} we obtain
\[\|f^e(x)+g^e(x)-h^e(x)\|\leq\varepsilon.\]
Thus we have
\[\|f^e(x+y)+g^e(x-y)-(f^e(x)+g^e(x))-(f^e(y)+g^e(y))\|\leq 3\varepsilon,\]
or equivalently,
\[\|f^e(x+y)+g^e(x-y)-f^e(x)-g^e(x)-f^e(y)-g^e(y)\|\leq 3\varepsilon.\]
Now let $u=f^e+g^e$ and $v=f^e-g^e$. Then
\begin{eqnarray*}
&&\|u(x+y)+u(x-y)-2u(x)-2u(y)\|\\&=&\|f^e(x+y)+g^e(x+y)+f^e(x-y)+g^e(x-y)\\
&&-2f^e(x)-2g^e(x)-2f^e(y)-2g^e(y)\|\\
&\leq& 6\varepsilon,
\end{eqnarray*}
and
\begin{eqnarray*}
\|v(x+y)-v(x-y)\|&=&\|f^e(x+y)-g^e(x+y)-f^e(x-y)+g^e(x-y)\|\\
&\leq&6\varepsilon,
\end{eqnarray*}
for all $x, y \in {\mathcal X}$ with $x\perp y$. Thus $u$ is an
approximately orthogonally quadratic mapping, $v$ is an
approximately orthogonally constant mapping, $f^e=\frac{u+v}{2}$
and $g^e=\frac{u-v}{2}$. Hence $f = f^o - \frac{1}{2}u -
\frac{1}{2}v$ and $g = g^o - \frac{1}{2}u + \frac{1}{2}v$.
\end{proof}

\bibliographystyle{amsplain}

\begin{thebibliography}{10}

\bibitem{ALO1} J. Alonso,
\textit{Some properties of Birkhoff and isosceles orthogonality
in normed linear spaces},
Inner product spaces and applications,
1-–11, Pitman Res. Notes Math. Ser., 376, Longman, Harlow, 1997.

\bibitem{ALO2} J. Alonso,
\textit{Uniqueness properties of isosceles orthogonality in normed
linear spaces},
Ann. Sci. Math. Qu´ebec 18 (1994), no. 1, 25–-38.

\bibitem{B-M1} C. Baak and M. S. Moslehian,
\textit{Stability of $J^*$-homomorphisms},
Nonlinear Analysis
(TMA), 63 (2005), 42--48.

\bibitem{B-M2} C. Baak and M. S. Moslehian,
\textit{On the stability of orthogonally cubic functional}, 
Kyungpook Math. J. (to appear). 


\bibitem{BIR} G. Birkhoff, {\it Orthogonality in linear metric
spaces}, Duke Math. J. {\bf 1} (1935), 169--172.

\bibitem {CZE} S. Czerwik,
{\it  Stability of Functional Equations of Ulam--Hyers--Rassias
Type}, Hadronic Press, Palm Harbor, Florida, 2003.

\bibitem{D-F-A} C. R. Diminnie, R. W. Freese and E. Z. Andalafte,
\textit{An extension of Pythagorean and isosceles orthogonality
and a characterization of inner-product spaces},
J. Approx. Theory
39 (1983), no. 4, 295–-298.

\bibitem {DRL} F. Drljevi\' c,
{\it On a functional which is quadratic on $A$-orthogonal
vectors}, Publ. Inst. Math. (Beograd) {\bf 54} (1986), 63--71.

\bibitem {FOC} M. Fochi,  {\it Functional equations in $A$-orthogonal vectors},
Aequationes Math. {\bf 38} (1989), 28--40.


\bibitem {HYE} D. H. Hyers,  {\it On the stability of the linear
functional equation},
 Proc. Nat'l. Acad. Sci. U.S.A. {\bf 27} (1941),
222--224.

\bibitem {H-I-R} D. H. Hyers, G. Isac and Th. M. Rassias,
 {\it Stability of Functional Equations in Several Variables},
 Birkh$\ddot{a}$user, Basel, 1998.

\bibitem{JAM} R. C. James,
{\it Orthogonality in normed linear spaces},
Duke Math. J. {\bf 12} (1945), 291--302.

\bibitem {JUN} S.-M. Jung,
{\it Hyers--Ulam--Rassias Stability of Functional Equations in
Mathematical Analysis},
Hadronic Press, Palm Harbor, Florida, 2001.

\bibitem{K-P} O. P. Kapoor and J. Prasad,
\textit{Orthogonality and characterizations of inner product
spaces},
Bull. Austral. Math. Soc. 19 (1978), no. 3, 403–-416.

\bibitem{M-M} M. Mirzavaziri and M. S. Moslehian,
\textit{A fixed point approach to stability of a quadratic equation},
 Bull. Braz. Math. Soc. (to appear). 

\bibitem{MOS1} M. S. Moslehian, {\it On the stability of the orthogonal Pexiderized
Cauchy equation},
J. Math. Anal. Appl., 318 (2006), no. 1, 221--223.

\bibitem{MOS2} M. S. Moslehian,
{\it On the orthogonal stability of the Pexiderized quadratic
equation},
J. Differ. Equations. Appl., {\bf 11} (2005), no. 11,
999--1004.

\bibitem {MOS3} M. S. Moslehian,
\textit{Approximately vanishing of topological cohomology groups},
J. Math. Anal. Appl., 318 (2006), no. 2, 758--771.

\bibitem {RAS1} Th. M. Rassias,
{\it  On  the stability of the linear
mapping in Banach spaces},
 Proc. Amer. Math. Soc. {\bf 72} (1978), 297--300.

\bibitem {RAS2} Th. M. Rassias(ed.), {\it Functional Equations,
Inequalities and Applications}, Kluwer Academic Publishers,
Dordrecht, Boston, Londonm, 2003.

\bibitem{SZA} Gy. Szab\' o, {\it Sesquilinear-orthogonally quadratic mappings},
Aequationes Math. {\bf 40} (1990), 190--200.

\bibitem {ULA} S. M. Ulam,
{\it Problems in Modern Mathematics},
Wiley, New York, 1960.

\bibitem{VAJ} F. Vajzovi\' c, {\it \" Uber das Funktional $H$ mit der Eigenschaft: $(x, y)=0
\Rightarrow H(x+y)+H(x-y)=2H(x)+2H(y)$}, Glasnik Mat. Ser. III
{\bf 22} (1967), 73--81.
\end{thebibliography}

\end{document}